\documentclass[12pt]{amsart}
\usepackage{amsbsy,amssymb,amscd,amsfonts,latexsym,amstext,delarray,
amsmath,epsfig,graphicx,bm, url} 
\input xypic
\usepackage{hyperref}

\newtheorem{thm}{Theorem}[section]
\newtheorem{prop}[thm]{Proposition}

\newtheorem{lem}[thm]{Lemma}
\newtheorem{defn}[thm]{Definition}
\newtheorem{rem}[thm]{Remark}

\numberwithin{equation}{section}

\def\bN{{\mathbb N}}

\def\Z{{\mathbb Z}}
\def\R{{\mathbb R}}

\def\fH{{\mathfrak H}}

\def\cutint{{\int \!\!\!\!\!\! -}}

\def\cA{{\mathcal A}}

\def\cL{{\mathcal L}}

\def\cO{{\mathcal O}}

\newcommand{\ie}{{\it i.e.\/}\ }

\newcommand{\cf}{{\it cf.\/}\ }

\def\text{\hbox}

\def\Aut{{\rm Aut}}

\def\Dom{{\rm Dom}}

\def\GL{{\rm GL}}

\def\Index{{\rm Index}}

\def\Tr{{\rm Tr}}

\def\vol{{\rm vol}}

\newcommand{\qqq}{{\,,\quad \forall\,}}

\def\Dirac{{\partial \hspace{-6pt} \slash}}
\begin{document}

\title {Type III and spectral triples}
\author[Connes]{Alain Connes}
\author[Moscovici]{Henri Moscovici}
\address{A.~Connes: Coll\`ege de France \\
3, rue d'Ulm \\ Paris, F-75005 France\\
I.H.E.S. and Vanderbilt University} \email{alain\@@connes.org}
\address{H.~Moscovici:
Department of mathematics, The Ohio State University, Columbus, OH 43210, USA}
\email{henri@math.ohio-state.edu}

\thanks{The work of the second named author was partially
 supported by the National Science Foundation
    award no. DMS-0245481}

\begin{abstract} We explain how a simple twisting of the notion of spectral
triple allows  to incorporate type III examples, such as those
arising from the transverse geometry of codimension one foliations.
We show that the classical cyclic cohomology valued Chern character
of finitely summable spectral triples extends to the twisted case
and lands in ordinary (untwisted) cyclic cohomology.  The index
pairing with ordinary (untwisted) K-theory continues to make sense
and the index formula is given by the pairing of the corresponding Chern
characters. This opens the road to extending the local index formula
to the type III case.
\end{abstract}

\maketitle

\section{Introduction}

The basic paradigm of noncommutative geometry is that of a spectral
triple $(\cA , \fH , D) $ (\cf \cite{book}, \cite{C-M2}), where the
algebra  $\cA$ encodes the space and the operator $D$ encodes the
metric. In the finite dimensional situation, \ie when there is an $\, \alpha > 0$ such that the  $n$-th
characteristic value of the resolvent of $D$ decays as
$n^{- \alpha}$ for $n\to \infty$, the Dixmier trace \cite{book} induces
a nontrivial trace on the algebra $\cA$. The existence of a trace is a
characteristic of the type II situation in the Murray-von Neumann
classification of rings of operators. Thus, in essence the theory is,
in its finite dimensional form, restricted to the type II case.
\smallskip

We shall explain in this note how a simple twisting of the notion of
spectral triple allows  to incorporate type III examples, such as
those arising from the transverse geometry of codimension one
foliations. Since the twisting of the commutators turns the usual hypertrace
constructed out of the Dixmier trace into a twisted trace on the coordinate algebra,
one would be tempted to interpret that as a manifestation of twisting
at the level of
cyclic cohomology, akin to that introduced by the authors in the context of Hopf
cyclic cohomology \cite{C-M}. The main point of this note, besides
giving simple natural examples of the general notion and developing
the first basic steps of the theory, is to show that contrary to
the initial expectations no cohomological twisting 
is in fact required. The Chern character of finitely summable spectral
triples extends to the twisted case, and lands in fact in ordinary
(untwisted) cyclic cohomology. The same holds true for the local
Hochschild character. The index
pairing with ordinary (untwisted) K-theory continues to make sense
and the index formula is still given by the pairing of the corresponding Chern
characters.
This opens the road to extending the
local index formula, as well as the analogue of the hypoelliptic construction on the dual
system together with the corresponding Thom isomorphism,  to the context of twisted
spectral triples of type III.
\medskip

\section{Two motivating examples}

\subsection{Dirac operator}
We start by recalling the classical comparison formula for the Dirac operators
associated to conformally equivalent metrics, \cf  \cite{hitchin}, \cite{bourg}.
Given a compact spin manifold $M^n$, to each Riemannian metric
 $g$ on $M$ one can canonically associate a Dirac operator $\Dirac = \Dirac^g$
 acting on the Hilbert space $\fH = \fH^g : = L^2(M,S^g)$ of $L^2$-sections of
 the spin bundle $S= S^g$,
and thus a corresponding spectral triple (cf. \cite{book})
$(\cA, \fH,\,  \Dirac)$ over the algebra $\cA: = C^\infty(M)$.
Let $h\in C^\infty(M)$ be a self-adjoint element and  replace $g$ by the
rescaled metric $g' = e^{-4h}\,g$. After identifying the corresponding spin bundles via
the $Spin_n$-equivariant transformation $\beta^g_{g'}$  from $g$-spinorial frames
 to $g^{\,\prime}$-spinorial frames defined in  \cite{bourg}, the gauge transformed operator
$\,^g \Dirac^{\,g'}\,  := \, \beta^{g'}_g \circ \Dirac^{\,g'} \circ
\beta^g_{g'}$ has the expression, \cf \cite[(26)]{bourg},
$$
^g \Dirac^{\,g'} \,  = \, e^{(n+1) h}\circ {\Dirac}^{\,g} \circ
e^{-(n-1)h} \, .
$$
In order to account for the change of the Riemannian volume form, $ \vol_{g'} = e^{-2nh}\,\vol_g$,
at the level of $L^2$-sections one needs to further rescale the canonical identification by
setting
$$
\tilde{\beta}^g_{g'} \, := \, e^{nh} \circ \beta^g_{g'} \, = \,
  \beta^g_{g'} \circ e^{nh}: \fH^g \rightarrow \fH^{g'} .
$$
This shows that the gauge transformed spectral triple is simply obtained by
replacing $\Dirac$ with
$$
{\Dirac}^{\,\prime} \, =\, e^{h}\,\Dirac \,e^{h} \, .
$$
\medskip

\subsection{Perturbed spectral triple} \label{example1}
In the general case, when
one starts from an arbitrary spectral triple  $(\cA , \fH , D) $
(\cf \cite{book}, \cite{C-M2}) and a
self-adjoint element $h = h^* \in \cA$, it is natural to wonder what are
the properties of the `perturbed' triple
\begin{equation}
(\cA , \fH , D') \,,\quad D'=\,e^{h}\,D\,e^{h} .\label{spectrip}
\end{equation}
The operator $D'$ is still self-adjoint but the basic boundedness
condition
\begin{equation}
[D,a] \quad \hbox{is bounded} \quad \forall \, a \in \cA \, ,
\label{boundedcom}
\end{equation}
will not necessarily hold, unless $h$ is in the center of $\cA$.

However, the following simple fact holds true.

\begin{lem}\label{conjugate}
For any self-adjoint element $h \in \cA$, letting
\begin{equation}
 \sigma(a)=\,e^{2h}\,a\,e^{-2h}, \quad  a \in \cA , \label{dfesigma}
\end{equation}
one has
\begin{equation}
d^{\,\prime}_\sigma a \, := \, D'\,a-\,\sigma(a)\,D' \quad \hbox{is
bounded} \quad \forall \, a \in \cA \, . \label{boundedcom1}
\end{equation}
\end{lem}

\proof The fact that the usual commutators
$[D,b]$, $b \in \cA$, are bounded implies the boundedness of
\begin{equation} \label{bound'}
d^{\,\prime}_\sigma a \, = \, e^{h}\,D\,e^{h}\,a\,
-\,e^{2h}\,a\,e^{-h}\,D\, e^{h} \, = \, e^{h} \, [D, b] \, e^{h} \,
, \qquad b=e^{h}\,a\,e^{-h} .
\end{equation}
\endproof

\medskip

\subsection{Transverse spectral triple } \label{trans1}
As a next example we take a codimension $1$ foliation and consider
the corresponding noncommutative algebra  $\cA$ of  `transverse coordinates'.
In fact, as a further simplification
we restrict to a complete transversal and take for $\cA$ the
algebraic crossed product of the algebra $C^\infty(S^1)$ of smooth
functions on $S^1$ by a group $\Gamma$ of orientation preserving
diffeomorphisms. Any element of $\cA=\,C^\infty(S^1)\rtimes \Gamma$
is represented as a finite sum of the form
$$ a=\,\sum_\Gamma
a_\phi\, U^*_\phi \,,
$$
the product rule is determined by
\begin{equation}
U^*_\phi\,f=\,(f \circ\phi)\, U^*_\phi \,,\quad
U^*_\phi\,U^*_\psi=U^*_{\psi\,\phi}\, , \label{product}
\end{equation}
and the involution
\begin{equation}
a=\,\sum_\Gamma
a_\phi\, U^*_\phi \quad \longmapsto \quad
a^* =  \,\sum_\Gamma U_\phi \, \bar{a}_\phi \, . \label{involution}
\end{equation}

One represents $\cA=\,C^\infty(S^1)\rtimes \Gamma$ in the Hilbert
space $\fH=\,L^2(S^1)$  by the
$\ast$-representation
\begin{equation}
(\pi
(g\,U^*_\phi)\,\xi)(x)=\,g(x)\,\phi'(x)^{\frac{1}{2}}\,\xi(\phi(x))\,
, \quad \forall \xi \in \fH \, , \quad x\in \R/\Z . \label{repa}
\end{equation}
In the role of $D$ we take the operator
$\, \Dirac =\,\frac{1}{i}\,\frac{d}{dx}$, while the automorphism
$\sigma \in \Aut \cA$ is defined on the monomials generating $\cA$ by
\begin{equation} \label{dfesigma1}
 \sigma(g\,U^*_\phi)=\,\frac{d\phi(x)}{dx}\,g\,U^*_\phi \, .
\end{equation}
One then has the following boundedness property.

\begin{lem}\label{foliation}
For any $\, a \in \cA$, the twisted commutators
\begin{eqnarray}
\Dirac \circ \pi (a) \,&-&\,\pi (\sigma(a))\circ \Dirac
  \, , \label{boundedcom2} \\
\text{\rm and} \qquad |\Dirac| \circ \pi (a) \,&-&\,\pi (\sigma(a))\circ |\Dirac|
 \label{boundedcom3}
\end{eqnarray}
are bounded.
\end{lem}

\proof For any $ a = g\,U^*_\phi$ one has:
\begin{eqnarray*}
\Dirac (\pi (a) \xi) (x)&= \frac{1}{i}\,\frac{d}{dx} \left(g(x)
\phi'(x)^{\frac{1}{2} }\right)\,\xi(\phi(x)) +
g (x) \phi'(x)^{\frac{1}{2} } \, \Dirac  (\xi \circ \phi)(x)\\
 &= \frac{1}{i}\,\frac{d}{dx} \left(g(x) \phi'(x)^{\frac{1}{2} }\right)\,\xi (\phi(x)) +
 g (x) \phi'(x)^{\frac{1}{2}} \, \left( \phi' (x) \, (\Dirac  (\xi) )(\phi(x)) \right) \, ;
\end{eqnarray*}
equivalently,
\begin{equation} \label{bdcom}
\Dirac  ((\pi (a)) (\xi) ) (x) -  \pi (\sigma (a)) (\Dirac
(\xi) ) (x) =  \frac{1}{i}\,\frac{d}{dx} \left(g(x) \phi'(x)^{\frac{1}{2}}\right)\,\phi'(x)^{-\frac{1}{2} }\, 
(\pi (U^*_\phi) \xi)(x),
 \end{equation}
which proves  the boundedness of the twisted commutator \eqref{boundedcom2}.
\smallskip

To prove the same for \eqref{boundedcom3} we shall switch from direct calculation
to an equally elementary symbolic argument.
Letting  $\, V_\phi$ denote the translation operator by  $\, \phi \in \Gamma$,
\begin{equation*}
(V_\phi\,\xi)(x)\, =\,\xi(\phi^{-1}(x))\,
, \quad \forall \xi \in \fH \, , \quad x\in \R/\Z \, ,
\end{equation*}
one has for any $ a = g\,U^*_\phi \in \cA$,
\begin{equation*}
\big( |\Dirac| \circ \pi(a)\, -\,\pi(\sigma(a)) \circ |\Dirac| \big)\circ V_\phi \, = \,
|\Dirac|  \circ g \,\phi'^{\frac{1}{2}}  \, - \, g \,\phi'^{\frac{1}{2}} \, \phi' \, V^{-1}_\phi \circ |\Dirac|  \circ V_\phi \, .
\end{equation*}
Now $V^{-1}_\phi \circ |\Dirac|  \circ V_\phi $ is a $1$-st order
(classsical) pseudodifferential operator whose principal symbol is
$\displaystyle \frac{1}{\phi'}$ times the principal symbol of  $|\Dirac| $.
It follows that the right hand
is a pseudodifferential operator of order $0$, hence bounded.
\endproof

  \begin{rem} \label{orderdecrease}
   {\rm The symbolic argument given above applies more generally
  to any pseudodifferential operator of arbitrary order $m \in \R$. Thus, if 
  $\, P \in \Psi DO^m (S^1)$, then for all $\, a \in \cA$}
 \begin{equation*}  \label{lowerorder}
\big( P \circ \pi(a)\, -\,\pi(\sigma^m(a)) \circ P \big)\circ V_\phi \ \in \, \Psi DO^{m-1} (S^1) \, .
\end{equation*}
  \end{rem}

\begin{rem} {\rm The canonical state $\varphi$ on $\cA$,
 \begin{equation}
\varphi(f \, U^*_\phi)=0 \;\textit{ if }\; \phi\neq 1\,,\quad \text{and} \quad
\varphi(f)=\int_{\R/\Z} \,f (x) \,dx \label{state}
\end{equation}
is a $\sigma^{-1}$-\textit{trace}, \ie satisfies 
\begin{equation} \label{sigmatrace}
\varphi(a\, b)\, =\, \varphi(b \, \sigma^{-1} (a))\, , \qquad \forall \, a, b \in \cA \, ,
\end{equation}
and its modular automorphism group  
is precisely the one-parameter group of automorphisms
\begin{equation}  \label{dfesigma1t}
 \sigma_t(g\,U^*_\phi)=\,\left(\frac{d\phi(x)}{dx}\right)^{it}\,g\,U^*_\phi\, , \qquad t \in \R \, ,
\end{equation}
whose value (after analytic continuation) at $\, t = -i$ coincides with
 $\sigma$.}
 \end{rem}

\medskip

\section{$\sigma$-spectral triples and their basic properties}

\subsection{Elementary properties} The usual definition of a spectral triple
extends to the context illustrated by the preceding examples as follows.

\begin{defn}
With $\sigma$ being an automorphism  of $\cA$, an {\rm ungraded} $\sigma$-spectral
triple $(\cA , \fH , D)$ is given by an action of
 $\cA$ in the Hilbert space $\fH$, while $D$ is a self-adjoint
operator with compact resolvent and such that
\begin{equation}
D\,a\, -\,\sigma(a)\,D \quad \hbox{is bounded} \quad \forall \, a \in \cA \, .
\label{boundedcom4}
\end{equation}
A {\rm graded} $\sigma$-spectral triple is similarly defined, with the additional
datum of a grading operator
$$ \gamma = \gamma^* \in \cL (\fH) \, , \quad \gamma^2 = I
$$
which commutes with the action of $\cA$, while
\begin{equation}
D\,\gamma \, = \, -\,\gamma\,D \, .
\label{anticomm}
\end{equation}
A  {\rm Lipschitz-regular} $\sigma$-spectral triple is one that
satisfies the additional condition
\begin{equation}
  \forall \, a \in \cA \, , \quad |D|\,a-\,\sigma(a)\,|D| \quad  \hbox{is bounded}.
\label{boundedcom5}
\end{equation}
\end{defn}
In the case when  $\cA $ is an involutive algebra and the
representation is involutive, to ensure 
 the compatibility between the
automorphism $\sigma$ and the $\ast$-involution, we impose the
additional \textit{unitarity condition}:
\begin{equation}
\sigma(a^*)=\,(\sigma^{-1}(a))^* \qqq \, a \in \cA \, .
\label{sigmainv}
\end{equation}

Lemma \ref{foliation} shows that the transverse Dirac operator on a
codimension $1$ foliation gives rise to a Lipschitz-regular
$\sigma$-spectral triple with $\sigma$ given by the Jacobian of the
holonomy. In particular, since this spectral triple is $1$-summable,
one gets examples of finitely summable $\sigma$-spectral triples for
which the representation of $\cA$ in $\fH$ generates a type III
factor.
\medskip

Let us spell out the extensions to the twisted case of some
basic properties of spectral triples. For background on spectral triples,
including notational conventions used below, we refer the reader
to \cite[IV.2]{book} and
\cite[Appendix A]{C-M2}, while for the notion of $\sigma$-trace see \cite{C-M}.

First of all, we note that any twisted spectral triple which is
Lipschitz-regular can be canonically `untwisted'  by passage to its
`phase'. This is quite clear in the case of the second
example (\cf  \S \ref{trans1}), since the phase of $\, \Dirac
=\,\frac{1}{i}\,\frac{d}{dx}$  is the Hilbert transform, and is actually easy
to prove in full generality.

\begin{prop} \label{untwist}
If the $\sigma$-spectral triple $(\cA , \fH , D)$ is
Lipschitz-regular and $\, F =  D \, |D|^{-1} $, then  $( \fH , F)$
is a Fredholm module over $\cA$. If moreover $(\cA , \fH , D)$ is
finitely summable, so is $( \fH , F)$.
\end{prop}

\proof
Indeed, for any $\, a \in \cA$,
\begin{equation*}
D\,a\, -\,\sigma(a)\,D \, = \, |D|\, (F\,a\, -\, a\, F ) \, + \, (|D|\,a-\,\sigma(a)\,|D| )\, F \, ,
\end{equation*}
therefore
\begin{equation*}
[F, a] \, = \, |D|^{-1} \, \big(  (D\,a\, -\,\sigma(a)\,D) \, -
\, (|D|\,a-\,\sigma(a)\,|D| )\, F \big) \, .
\end{equation*}
Thus, all these commutators are compact operators, and in fact they are
quantized differentials of the same order as $\, D^{-1}$.
\endproof
\medskip

\begin{prop} \label{resandsig}
Let $(\cA , \fH , D)$ be a $\sigma$-spectral triple with
$D^{-1}\in \cL^{n,\infty}$.
\begin{itemize}
\item[(1)] The linear functional
\begin{equation*}
a \in \cA \mapsto \varphi(a)=\,\cutint \,a\,D^{-n} \, := \Tr_\omega (a\,D^{-n})
\end{equation*}
is a $\sigma^{-n}$-trace on $\cA$:
$\quad \varphi(a\, b)\, =\, \varphi(b \, \sigma^{-n} (a))$, $\quad \forall \, a, b \in \cA$.
\item[(2)] More generally, for any bounded operator
$T \in \cL (\fH)$,
\begin{equation}\label{dixsig2}
Tr_\omega (T \, \sigma^{-n} (a) \,D^{-n}) \, = \, Tr_\omega (a\, T \, D^{-n} ) \, ,
\qquad \forall \, a \in \cA \, .
\end{equation}
\item[(3)] When the $\sigma$-spectral triple $(\cA , \fH , D)$ is Lipschitz-regular,
the same hold true when $\,D^{-n}$ is replaced by $\,|D|^{-n}$.
\end{itemize}
\end{prop}

\proof Let us show by induction that for any $\, 1 \leq k \leq n$ one has
\begin{equation} \label{ord-k}
D^{-k}\,a-\,\sigma^{-k}(a)\,D^{-k}\in \cL_0^{\frac{n}{k} \, ,\infty}\qqq a \in \cA .
\end{equation}
Clearly,
$$
D^{-1}\,a-\,\sigma^{-1}(a)\,D^{-1}=\,D^{-1}\,(a\,D-\,D\,\sigma^{-1}(a))\,D^{-1}  \, \in  \, \cL_0^{n ,\infty} .
$$
To verify the inductive step we write
\begin{eqnarray*}
D^{-k}\,a-\,\sigma^{-k}(a)\,D^{-k} &=& D^{-1} \big(D^{-(k-1)}\,a-\,\sigma^{-(k-1)}(a)\,D^{-(k-1)}\big)\\
&+& \big(D^{-1}\,\sigma^{-(k-1)}(a)\, -\,\sigma^{-k}(a)\,D^{-1}\big) D^{-(k-1)}
\end{eqnarray*}
and observe that, by H\"older's inequality  and the induction hypothesis, each of the two
summands in the right hand side belongs to $\,\cL_0^{\frac{n}{k} \, ,\infty}$.

Applying now \eqref{ord-k} for $k=n$, one obtains
$$
\varphi(T\,\sigma^{-n}(a))=\,\cutint \,T\,\sigma^{-n}(a)\,D^{-n}=\,
\cutint \,T\,D^{-n}\,a=\, \varphi(a\,T)  \, , \quad \forall \, T \in \cL (\fH) \, .
$$
To prove the third statement, one replaces $\,D$ by $\,|D| $ throughout the above argument.
\endproof

We now consider the analogue of the bimodule  of gauge potentials
given in the usual case by  the $\cA$-bimodule $\Omega_D^1 \subset
\cL (\fH)$
of operators of the form
\begin{equation}\label{pot}
A = \Sigma  a_i [D,b_i], \ \ \  a_i , b_i \in \cA \, .
\end{equation}
Let $(\cA , \fH , D)$ be a $\sigma$-spectral triple, then we let
$\Omega_D^1 \subset \cL (\fH)$ be the linear space of operators of
the form
\begin{equation}\label{pot1}
A = \Sigma  a_i (D\,b_i-\sigma(b_i)\,D)\, , \ \ \  a_i , b_i \in \cA
\, .
\end{equation}
\medskip

\begin{prop} \label{bim}
Let $(\cA , \fH , D)$ be a $\sigma$-spectral triple, then
$\Omega_D^1 $ is an $\cA$-bimodule for the action
\begin{equation}\label{bim1}
a \cdot\omega\cdot b=\,\sigma(a)\,\omega\,b \qqq a,\,b \in \cA \qqq
\omega \in \Omega_D^1
\end{equation}
and the map
\begin{equation}\label{der}
a \mapsto d_\sigma (a)=\,D\,a-\sigma(a)\,D
\end{equation}
is a derivation of $\cA$ in $\Omega_D^1 $.
\end{prop}

\proof One has
$$
d_\sigma (ab)=\,D\,ab-\sigma(ab)\,D=\,(D\,a-\sigma(a)\,D)\,b+\sigma(a)\,(D\,b-\sigma(b)\,D)
$$
which shows that
\begin{equation}\label{bim2}
d_\sigma (ab)=\,d_\sigma (a)\cdot b+ a\cdot d_\sigma (b) \qqq a, b \in \cA
\end{equation}
Since $\sigma$ is an automorphism of $\cA$, the linear space
$\Omega_D^1$ is the linear span of the $a\cdot d_\sigma (b)$ for $a, b \in
\cA$. By \eqref{bim2} this is stable under right multiplication by
elements of $\cA$. Thus $\Omega_D^1$ is an $\cA$-bimodule. Finally
\eqref{bim2} shows that $d_\sigma $ is a derivation.
\endproof

\medskip

\subsection{Chern character} By Proposition \ref{untwist},
any $\sigma$-spectral triple $(\cA , \fH , D)$ of finite summability
degree, \ie such that with $D^{-1}\in \cL^{n,\infty}$ for some $n
\in \bN$, which in addition is {\it Lipschitz-regular} has a
well-defined Chern character in cyclic cohomology, namely the Chern
character of its `phase' Fredholm module $( \fH , F)$ over $\, \cA$,
\begin{equation}\label{CCh}
\Phi_F (a^0 , a^1, \ldots , a^n)\, :=  \,
\Tr \, (\gamma\, F\,  [F, a^0] \, [F, a^1] \cdots \, [F, a^n]) \, ,
\qquad  \forall \, a^0 , a^1, \ldots a^n \in \cA ,
\end{equation}
with $\gamma$ omitted in the ungraded case (\cf~\cite[Part I]{ncdg}).

On the other hand
let us assume for a moment
that $(\cA , \fH , D')$, $D' = e^h \, D \, e^h$,
is a graded twisted spectral triple as in \S \ref{example1}, with the
property that $D^{-1}\in \cL^{n,\infty}$ for some even $n \in \bN$.
Applying \eqref{bound'} one sees that
\begin{equation*}
D'^{-1} \, d'_\sigma a \, = \, e^{-h} \, D^{-1} \, [D, b] \, e^h \, ,  \quad \text{where} \quad
b = e^h \, a \, e^{-h} .
\end{equation*}
Therefore, for any $a^0 ,  a^1, \ldots ,  a^n \in \cA$,
\begin{equation} \label{compchern}
\Tr \, (\gamma\, D'^{-1}  d'_\sigma a^0 \, D'^{-1} d'_\sigma a^1 \cdots \, D'^{-1} d'_\sigma a^n)
=  \Tr \, (\gamma\, D^{-1} [D, b^0] \, D^{-1} [D, b^1] \cdots \,
D^{-1} [D, b^n]) \, ,
\end{equation}
with $\, b^i \, =\, e^h \, a^i \, e^{-h}$, $ \forall \, i = 0, \ldots , n$.
The right hand side of the above identity is a cyclic cocycle on $\cA$ that
represents, up to normalization, the Chern character of
$(\cA , \fH , D)$, \cf~\cite[Part I, \S 6]{ncdg}. It follows that the left hand side
is also a cyclic cocycle, obtained via conjugation by an inner automorphism,
and thus determining the same periodic cyclic cohomological class
\begin{equation*}
 Ch^* (\cA , \fH , D') \, = \,  Ch^* (\cA , \fH , D) \in HP^* (\cA) .
\end{equation*}

This suggests that it should be possible to define a `straight'
Chern character for any finitely sumable twisted spectral triple,
not just for those that are Lipschitz-regular. The proposition below
confirms that this is indeed the case.
\medskip

\begin{prop} \label{straightchern}
Let $(\cA , \fH , D)$ be a graded $\sigma$-spectral triple such that
$D^{-1}\in \cL^{n,\infty}$ for some even $n \in \bN$. Then the following multilinear form
\begin{equation}\label{Cycocycle}
\Phi_{D, \sigma} (a^0 , a^1, \ldots , a^n)\, :=  \,
\Tr \, (\gamma\, D^{-1}  d_\sigma a^0 \, D^{-1} d_\sigma a^1 \cdots \, D^{-1} d_\sigma a^n),
\quad  a^0, \ldots a^n \in \cA
\end{equation}
is a cyclic cocycle in $Z_\lambda^n (\cA)$.
\end{prop}

\proof The proof of Proposition 1 in~\cite[Part I, \S 6]{ncdg} (\cf also  {\it infra})
applies verbatim
to the twisted case if one simply replaces the representation
$\, a \mapsto D^{-1} \, a \, D\, $ by  $\, a \mapsto D^{-1} \, \sigma(a) \, D$,
$\, \forall \, a \in \cA$.
\endproof

\medskip

\subsection{Pairing with $K$-theory} Implicit in the above proof is the
existence `behind the scene' of a pair of graded Fredholm module over $\cA$,
$\, (\tilde\fH^{\pm}, F^{\pm})\,$,
canonically associated
to the $\sigma$-spectral triple $(\cA , \fH , D)$. To wit, we
decompose $(\cA , \fH , D)$
according to the grading by $\gamma$ into
\begin{equation*}
\fH =  \fH_{+} \oplus \fH_{-} \, , \qquad
D  = \begin{pmatrix} 0  &D_{-}\\ \, D_{+}   &0  \end{pmatrix} \, , \qquad
a  =  \begin{pmatrix}  \, a_{+}  &0\\  0   &a_{-}  \end{pmatrix} \, ,
\quad \forall a \in \cA ,
\end{equation*}
then define
\begin{equation*}
\tilde\fH^{\pm}  :=  \fH_{\pm}  \oplus \fH_{\pm} , \quad
\pi^{\pm} (a)  :=
\begin{pmatrix}   a_{\pm}   &0\\ 0 &D_{\pm}^{-1}  \sigma(a_{\mp})  D_{\pm}  \end{pmatrix} \, \,
\text{on}\, \, \Dom D_{\pm} , \quad
F^{\pm} := \, \begin{pmatrix} 0  &I_{\pm}\\ \, I_{\pm}   &0   \end{pmatrix},
\end{equation*}
and note that for all $\,   a \in \cA \,$, firstly,
\begin{equation} \label{Fred1}
D_{\pm}^{-1}  \sigma(a_{\mp})  D_{\pm}  \,= \, a_{\pm} \,
- \, D_{\pm}^{-1} (D_{\pm}\, a_{\pm} \, - \, \sigma (a_{\mp}) \, D_{\pm})  \quad
\text{is bounded} ,
  \end{equation}
and secondly,
\begin{equation} \label{Fred2}
[F^{\pm} , \pi^{\pm} (a)] \, =
\begin{pmatrix}   0  &- D_{\pm}^{-1} \left(D_{\pm} \, a_{\pm} - \sigma(a_{\mp}) \, D_{\pm}\right) \\
D_{\pm}^{-1} \left(D_{\pm} \, a_{\pm} - \sigma(a_{\mp}) \, D_{\pm}\right) &0\end{pmatrix}
\in \cL^{n,\infty}.
 \end{equation}
\medskip

\begin{lem} \label{pairing}
Let $(\cA , \fH , D)$ be a graded $\sigma$-spectral triple such that
$D^{-1}\in \cL^{n,\infty}$ for some even $n \in \bN$, and let
$e \in  \cA$ be an  idempotent. Denote by $f_{\pm}$ the bounded closure of
of $D_{\pm}^{-1}  \sigma(e_{\mp})  D_{\pm}$. Then
$\, f_{\pm}^2 = f_{\pm}\,$  and
$\, f_{\pm} e_{\pm}  : e_{\pm} \fH_{\pm} \longrightarrow f_{\pm} \fH_{\pm}\, $
are Fredholm operators.
\end{lem}

\proof The first claim is obvious and the second follows from the fact that
$\, f_{\pm}\, - \, e_{\pm}\,$ is compact.
\endproof
\medskip

The integer $\,\Index (f_{\pm} e_{\pm} )$ depends only on the $K$-theory
class of the idempotent, and  thus one can define a pair of index maps
$\, \Index^{\pm}_{D, \sigma} : K_0 (\cA) \rightarrow \Z $, by setting
\begin{equation} \label{indmap}
\Index^{\pm}_{D, \sigma} [e] \, =  \, \Index (f_{\pm} e_{\pm} ) \, , \qquad
\forall \, e^2 = e \in M_N (\cA) ;
\end{equation}
taken together, they give rise to a {\it double index map}
\begin{equation} \label{doublemap}
 \Index_{D, \sigma} = \big(\Index^{+}_{D, \sigma} ,  \, \Index^{-}_{D, \sigma} \big)
: K_0 (\cA) \rightarrow \Z \times \Z .
\end{equation}
On the other hand, the cyclic cocycle \eqref{Cycocycle} is itself made
of two cocycles in $Z_\lambda^n (\cA)$,
\begin{equation}\label{pmcocycle}
\Phi^{\pm}_{D, \sigma} (a^0, \ldots , a^n) := 
\Tr \,\big( D_{\pm}^{-1} (D_{\pm} \, a^0_{\pm} - \sigma(a^0_{\mp}) \, D_{\pm})\,
 \cdots \, D_{\pm}^{-1} (D_{\pm} \, a^n_{\pm} - \sigma(a^n_{\mp}) \, D_{\pm})\big).
\end{equation}
\medskip

\begin{prop} \label{doubleindex}
Let $(\cA , \fH , D)$ be a graded $\sigma$-spectral triple such that
$D^{-1}\in \cL^{n,\infty}$ for some even $n \in \bN$.
For any $\, e^2 = e \in M_N (\cA) $, one has
\begin{equation}\label{doubleindexformula}
\Index^{\pm}_{D, \sigma} [e] \, =  \, \Phi^{\pm}_{D, \sigma} (e , \ldots , e) \, .
\end{equation}
If in addition $\, e^* = \sigma (e) $, then
\begin{equation} \label{sigidemp}
\Index^{+}_{D, \sigma} [e] \, =  \, - \, \Index^{-}_{D, \sigma} [e] \, .
\end{equation}
\end{prop}

\proof
Indeed, $\, \Index^{\pm}_{D, \sigma} $ is precisely the index map associated
to the Fredholm module $\, (\tilde\fH_{\pm} , F_{\pm} )\,$, and therefore
is given by the corresponding index formula, \cf~\cite[Part I, \S 3, Theorem 1]{ncdg}.
\smallskip

The second claim follows from the fact that, if  $\, e^* = \sigma (e)$, then
\begin{equation}  \label{sigadj1}
\big(D  \, e - \sigma(e) \, D\big)^* \, = \, - \, (D  \, e - \sigma(e) \, D) \, ,
\end{equation}
which in turn implies
\begin{equation*}
\overline{\Phi^{+}_{D, \sigma} (e , \ldots , e)} \, =  \, - \Phi^{-}_{D, \sigma} (e , \ldots , e)  \, .
\end{equation*}
\endproof

\medskip

\subsection{One-parameter group of automorphisms} We now make the the additional
assumption that
the involutive Banach algebra $\cA$ is equipped with a strongly continuous
$1$-parameter group of isometric automorphisms $\,  \{ \sigma_t \}_{t \in \R} $
such that
\begin{itemize}
\item[(1PG)]  {\em $\sigma$ coincides with the value at $t = -i$ of the analytic extension
of $ \{ \sigma_t \}_{t \in \R}$}.
\end{itemize}
The existence of such an analytic extension defined on a {\em dense subalgebra}
$\cO$ of $ \cA$, which is moreover {\em stable under holomorphic functional calculus}, is
ensured by a theorem of Bost~\cite[Thm. 1.1.1]{bost}.
\smallskip

To begin with, let us note that in the presence of the above hypothesis, which by the way is
automatically satisfied by the transverse spectral triple, \cf  \eqref{dfesigma1t},
the double index map reduces to a single index pairing map.
\medskip

\begin{lem} Assume that $\cA$ satisfies {\rm (1PG)}.
Then the two signed index maps coincide, \ie
\begin{equation*}
\Index^{+}_{D, \sigma} \, =  \, - \, \Index^{-}_{D, \sigma} : \, K_0 (\cA) \rightarrow \Z  \, .
\end{equation*}
\end{lem}

\proof Let $\, e^2 = e = e^* \in M_N (\cO) $ be a projection. Instead of \eqref{sigadj1}
we now use the identity
\begin{equation}  \label{sigadj2}
\big(D  \, \sigma^{-1}(e)\, - \, e \, D\big)^* \, = \, - \big( D  \, e \,  - \,\sigma (e) \, D \big) \, ,
\end{equation}
together with the fact that the idempotents $e$ and $\sigma^{-1} (e)$ are homotopic
via the path $ t \in [0, 1] \mapsto \sigma_{it} (e)$, to obtain
\begin{eqnarray*}
\Index^{+}_{D, \sigma} [e] &=&
\Index^{+}_{D, \sigma} [\sigma^{-1}(e)] =
\overline{\Phi^{+}_{D, \sigma} (\sigma^{-1}(e) , \ldots , \sigma^{-1}(e))} =
- \Phi^{-}_{D, \sigma} (e , \ldots , e) \\
 &=& - \Index^{-}_{D, \sigma} [e]  \, .
\end{eqnarray*}
\endproof

As a matter of fact, an elaboration of the above argument gives a more
comprehensive result that elucidates the relationship between the two
`half-character' cocycles \eqref{pmcocycle}. Note however that
unlike the classical case, where this relationship manifests itself already at the level of
cocycles (\cf~\cite[Part I, \S 6]{ncdg}), in the twisted case it only occurs at the
cohomological level.
\medskip

\begin{thm}
Let $(\cA , \fH , D)$ be a graded $\sigma$-spectral triple such that
$D^{-1}\in \cL^{n,\infty}$ for some even $n \in \bN$,
Under the assumption {\rm (1PG)}, the two  Chern characters
$[\Phi^{\pm}_{D, \sigma}]  \in HP^{\rm ev} (\cO)$
are related by the identity
\begin{equation} \label{adjointchern}
[\Phi^{-}_{D, \sigma}]\,  =\, - [(\Phi^{+}_{D, \sigma})^*]  \, ,
\end{equation}
where
\begin{equation*}
(\Phi^{+}_{D, \sigma})^* (a_0 ,  \ldots , a_n) \, := \,
\overline{\Phi^{+}_{D, \sigma} (a_n^* , \ldots , a_0^*)} \, , \qquad \forall \,
a^0 , a^1, \ldots , a^n \in \cO \, .
\end{equation*}
\end{thm}

\proof
For any $\, a \in \cO$ and $\, t \in [0, 1]$,  let
\begin{equation} \label{homotop}
\pi_t  (a) \, := \, \big((\pi^{+} \circ \sigma_{it}) (a^*)\big)^* \,= \,
\begin{pmatrix}   \sigma_{-it} (a_{+})   &0\\ 0 &D_{-} \sigma_{i-it} (a_{-})  D_{-}^{-1}   \end{pmatrix} ,
\end{equation}
and define the family of Fredholm modules
$\{ (\fH_t , F_t )\}_{t \in [0, 1]}$, by taking
$$
\fH_t  =  \fH_{+}  \oplus \fH_{+}  \quad \text{acted upon by} \,\,  \cO \, \text{via} \,\,
\pi_t \, , \, \, \text{and} \quad
F_t := \, \begin{pmatrix} 0  &I_{+}\\ \, I_{+}   &0   \end{pmatrix}.
$$
As in \eqref{Fred2} one has
$\, [F_t , \pi_t (a)] \, \in \cL^{n,\infty}$, because
\begin{equation*}
 D \, \sigma_{i-it} (a) \, D^{-1} \, - \,  \sigma_{-it} (a)  \, = \,
\big( D \, \sigma_{i-it} (a) \, - \,  \sigma(\sigma_{i-it} (a) ) \, D\big)\, D^{-1}
\, \in \cL^{n,\infty}.
 \end{equation*}
Note also that
\begin{equation}   \label{pi0}
\pi_0 (a) \,= \,
\begin{pmatrix}  a_{+}   &0\\ 0 &D_{-} \sigma^{-1}(a_{-})  D_{-}^{-1}   \end{pmatrix}
  \,= \, \big(\pi^{+} (a^*)\big)^* \, ,
 \end{equation}
and
\begin{equation}  \label{pi1}
\pi_1 (a) \,= \,
\begin{pmatrix}   \sigma (a_{+})   &0\\ 0 &D_{-} a_{-}  D_{-}^{-1}   \end{pmatrix}
  \,= \,  \begin{pmatrix} 0  &D_{-} \\ \, D_{-}   &0  \end{pmatrix} \, \pi^{-} (a)
  \,  \begin{pmatrix} 0  &D_{-}^{-1} \\ \, D_{-}^{-1}   &0  \end{pmatrix} .
\end{equation}
By Lemma 1 in~\cite[Part I, \S 5]{ncdg}, the periodic cyclic cohomology class
$\, Ch^* (\fH_t , F_t ) \in HP^{\rm ev} (\cO)$,  is independent of $t \in [0, 1]$.
We also recall that this class can be represented by the cyclic cocycle
\begin{equation*}
\Phi_t (a^0 ,  \ldots , a^n) := 
\frac{1}{2} \, \Tr \left(\begin{pmatrix} I  &\, 0\\ \, 0   &-I   \end{pmatrix}
F_t \, [F_t , \pi_t (a^0)] \, \cdots \, [F_t, \pi_t (a^n)] \right), \quad
a^0 , \ldots , a^n \in \cO .
\end{equation*}
In particular, using \eqref{pi0} and \eqref{sigadj2}, one has
\begin{eqnarray*}
\Phi_0 (a^0 ,  \ldots , a^n )&=& \,
 \Tr \left( (a^0_{+} - D_{-} \sigma^{-1}(a^0_{-})  D_{-}^{-1} ) \, \cdots \,
(a^n_{+} - D_{-} \sigma^{-1}(a^n_{-})  D_{-}^{-1} ) \right) \\
&=&  \,  (\Phi^{+}_{D, \sigma})^*  (a^0 ,  \ldots , a^n ) \, ,
\end{eqnarray*}
while by \eqref{pi1}
\begin{eqnarray*}
\Phi_1 (a^0 ,  \ldots , a^n ) \, &=& \,
 \Tr \left( (\sigma (a^0_{+} )- D_{-} a^0_{-}  D_{-}^{-1} ) \, \cdots \,
(\sigma (a^n_{+} )- D_{-} a^n_{-}  D_{-}^{-1} ) \right) \\
&=&\,   - \Phi^{-}_{D, \sigma}  (a^0 ,  \ldots , a^n ) .
\end{eqnarray*}
\endproof

\medskip

\subsection{Local Hochschild cocycles} \label{hoch}
As a first step in the direction of extending
the local index formula of~\cite{C-M2} to twisted spectral triple, we shall construct
an analogue of the local Hochschild cocycle that
gives the Hochschild class of the Chern character in the untwisted case,
\cf~\cite[IV.2.$\gamma$]{book}.
We begin by revisiting the latter and then proceed in a heuristic manner.

Given a graded spectral triple $(\cA , \fH , D)$ such that
$D^{-1}\in \cL^{n,\infty}$ for some $n \in 2 \bN$, the Hochschild class of
its Chern character is represented in local form by the following cocycle:
\begin{equation}\label{Hcocycle1}
\Psi_D (a^0 , a^1, \ldots , a^n)\, :=
\,\cutint \, \gamma \,a^0 \, [D, a^1]\, \cdots \, [D, a^n] \,D^{-n} \, ,
\qquad  \forall \, a^0 , a^1, \ldots , a^n \in \cA \, .
\end{equation}
Noting that
\begin{equation}\label{switch1}
[D, a] \,D^{-k} \, = \, D^{-k+1} \, \big(D^k \, a \, D^{-k} \, - \, D^{k-1} \, a \, D^{-k +1}\big) \, ,
\qquad  \forall \, a \in \cA \, ,
\end{equation}
we can successively move $\,D^{-n} \,$ to the left with loss of
a power at each step, and thus rewrite the cocycle \eqref{Hcocycle1}
in the form
\begin{equation}\label{Hcocycle2}
\Psi_D (a^0 , a^1, \ldots , a^n) \, =  \,
 \cutint \gamma \, a^0 \, (D \, a^1\, D^{-1} \, - \,  a^1) \, \cdots \,
(D^n \, a^n \, D^{-n} \, - \, D^{n-1} \, a^n\, D^{-n +1}) .
\end{equation}
In the twisted case, inspired by the formulas \eqref{pi0}, \eqref{pi1},  we make
the formal substitution
\begin{equation}\label{subst}
D^k \, a \, D^{-k}  \quad \longmapsto \quad
D^k \, \sigma^{-k} (a) \, D^{-k}  \, , \qquad  \forall \, a \in \cA \, ,
\end{equation}
and obtain the following candidate for a Hochschild character cocycle:
\begin{displaymath}
\begin{split}
\Psi_{D, \sigma}& (a^0 , a^1, \ldots , a^n) \, =  \\
 &\cutint \gamma \, a^0 \, \big(D \, \sigma^{-1} (a^1)\, D^{-1} \, - \,  a^1\big) \, \cdots \,
\big(D^n \, \sigma^{-n} (a^n) \, D^{-n} \, - \, D^{n-1} \, \sigma^{-n+1} (a^n)\, D^{-n +1}\big) .
\end{split}
\end{displaymath}
The counterpart of \eqref{switch1} being
\begin{equation}\label{switch2}
d_\sigma (\sigma^{-k} (a)) \,D^{-k} \, =
\, D^{-k+1} \, \big(D^k \, \sigma^{-k} (a) \, D^{-k} \, - \, D^{k-1} \, \sigma^{-k+1} (a) \, D^{-k +1}\big) \, ,
\end{equation}
one can reverse the process of distributing $\,D^{-n} \,$ among the factors, which
leads to the expression stated below.
\medskip

\begin{prop} \label{tHclass}
Let $(\cA , \fH , D)$ be a graded $\sigma$-spectral triple such that
$D^{-1}\in \cL^{n,\infty}$ for some even $n \in \bN$. Then the $n+1$-linear form
on $\cA$
\begin{equation}\label{tHcocycle}
\Psi_{D, \sigma} (a^0 , a^1, \ldots , a^n)\, :=
\,\cutint \, \gamma \,a^0 \, d_\sigma (\sigma^{-1}(a^1))\, \cdots \,
d_\sigma (\sigma^{-n} (a^n)) \,D^{-n} 
\end{equation}
is a Hochschild cocycle in $Z^n (\cA, \cA^*)$.

In the ungraded case, for  a $\sigma$-spectral triple of odd summability degree
which is Lipschitz-regular, the corresponding
Hochschild cocycle is defined by the expression
\begin{equation}\label{oddtHcocycle}
\Psi_{D, \sigma} (a^0 , a^1, \ldots , a^n)\, :=
\,\cutint \,a^0 \, d_\sigma (\sigma^{-1}(a^1))\, \cdots \,
d_\sigma (\sigma^{-n} (a^n)) \,|D|^{-n} \, .
\end{equation}
\end{prop}

\proof  We check that $\Psi_{D, \sigma}$ is a Hochschild cocycle by computing
its coboundary,
using the derivation rule \eqref{bim2}, as follows:
\begin{displaymath}
\begin{split}
&b\Psi (a^0 , a^1, ..., a^{n+1})= \sum_{i=0}^n (-1)^i\, \Psi (a^0,..., a^i a^{i+1},..., a^{n+1})
+ (-1)^{n+1} \Psi (a^{n+1} a^0 , a^1, ..., a^n) \\
 &=\, \cutint \, \gamma\, a^0 \, a^1\, d_\sigma (\sigma^{-1} (a^2)) \cdots \,
 d_\sigma (\sigma^{-n}  (a^{n+1}))\,D^{-n} \,\\
 &\quad - \,  \cutint \, \gamma \, a^0 a^1 \,
 d_\sigma (\sigma^{-1} (a^2)) \cdots \,   d_\sigma  (\sigma^{-n}  (a^{n+1}))\,D^{-n} \\
&\quad  - \, \cutint \, \gamma \, a^0\, d_\sigma (\sigma^{-1} (a^1)) \, \sigma^{-1} (a^2) \cdots \,
 d_\sigma (\sigma^{-n}  (a^{n+1}))\,D^{-n} \, + \ldots \\
&\quad \ldots  + \,   (-1)^n \cutint \, \gamma\, a^0\, d_\sigma (\sigma^{-1} (a^1)) \cdots
 \, \sigma^{-n+1} (a^{n-1})\, \sigma^{-n+1}  (a^n) \, \sigma^{-n} (a^{n+1}) \,D^{-n} \\
&\qquad \qquad + \,   (-1)^n \cutint \, \gamma \,a^0\, d_\sigma (\sigma^{-1} (a^1)) \cdots
 \, \sigma^{-n+1} (a^{n-1})\, d_\sigma ( \sigma^{-n}(a^n) )\,  \sigma^{-n}(a^{n+1}) \,D^{-n} \\
&\qquad  \qquad +\, (-1)^{n+1} \cutint  \, \gamma \,  a^{n+1} \, a^0\,
d_\sigma (\sigma^{-1} (a^1))  \cdots \, d_\sigma ( \sigma^{-n}(a^n) )\,D^{-n} \, .
\end{split}
\end{displaymath}
The resulting expression vanishes because of successive cancelations, with the last two terms
canceling each other in view of \eqref{dixsig2}.

In the ungraded Lipschitz-regular case, the very same calculation
holds true provided $\gamma$ is replaced by the phase operator $\, F
= \, D \, |D|^{-1}$.
\endproof
\medskip

 The cyclic group generated by $\sigma \in \Aut (\cA)$ acts in a
 natural way on the set of such Hochschild cocycles.
 For each integer $m \in \Z$, the corresponding `gauge transformed'
 cocycle via the action of $\sigma^m \in \Aut (\cA)$,  has the expression:
$\quad  \forall \, a^0 , a^1, \ldots a^n \in \cA$,
\begin{equation}\label{mtHcocycle}
\Psi^{(m)}_{D, \sigma} (a^0 , a^1, \ldots , a^n)\, :=
\,\cutint \, \gamma \,\sigma^{m}(a^0) \, d_\sigma (\sigma^{m-1}(a^1))\, \cdots \,
d_\sigma (\sigma^{m-n} (a^n)) \,D^{-n} \, .
 \end{equation}
 
\medskip

As an illustration, let us specialize formula \eqref{oddtHcocycle}
to the twisted spectral triple $(\cA, \Dirac, \fH)$
associated to a codimension 1 foliation (\cf \S \ref{trans1}). We shall assume
throughout  that
the group $\Gamma$ consists of orientation preserving
diffeomorphisms with \emph{nondegenerate isolated fixed points}.
Denoting the local Hochschild cocycle by $\Psi_1$, and using \eqref{bdcom}, one has 
\begin{eqnarray*}
&&\Psi_1 (f\,U^*_\phi , \, g\, U^*_\psi )  \, = \,
\cutint \,\pi (f\, U^*_\phi)
 \, \big( \Dirac \, \pi(\sigma^{-1}(g\, U^*_\psi)) \, - 
 \, \pi (g\, U^*_\psi) \,  \Dirac \big) \,|\Dirac|^{-1} \\ 
 &&= \,\cutint \, f  \, (\phi')^{\frac{1}{2}}\,  V^{-1}_\phi \, \frac{1}{i} \frac{d}{dx}
 \big(g\, (\psi')^{-\frac{1}{2}} \big) \,  V^{-1}_\psi \,|\Dirac|^{-1} 
  = \frac{1}{i} \,\cutint \, f  \, (\phi')^{\frac{1}{2}}\,
 \big(g\, (\psi')^{-\frac{1}{2}} \big)^\prime \circ \phi\, \, V^{-1}_{\psi\circ \phi} \,|\Dirac|^{-1} .
  \end{eqnarray*}
Here the functional $\displaystyle \cutint$ is extended to the algebra 
$\Psi DO^{\infty} (S^1) \rtimes \Gamma$  by the formula
\begin{equation} \label{res}
\cutint \, V^{-1}_\chi\, P \, = \, {\rm Res}_{s=0} \, \Tr (V^{-1}_\chi\, P \, |\Dirac|^{-s} )  \, ,
\qquad P \in \Psi DO^{\infty} (S^1) \, , \quad \chi \in \Gamma.
  \end{equation}
By~\cite[IV. 2.$\beta$]{book}, when  $\chi = I$
this coincides with the Wodzicki residue of  $P$. On the other hand,
when $P = F \, |\Dirac|^{-1} $, with $ F \in C^\infty (S^1)$, and $\chi$  
is a diffeomorphism with nondegenerate isolated fixed points
the above expression vanishes. 
Indeed, using the Mellin transform 
\begin{equation*} 
  \Tr (V^{-1}_\chi\, F \, |\Dirac|^{-s-1} )  \, = \, \Tr \left(V^{-1}_\chi\, F \, 
  \left(\Dirac^2\right)^{-\frac{s+1}{2}} \right)   \, = \,
  \frac{1}{\Gamma (\frac{s+1}{2} )} \, \int_0^\infty \, t^{\frac{s-1}{2}} \, 
  \Tr \left(V^{-1}_\chi\, F \, e ^{-t \, \Dirac^2} \right) \, dt .
   \end{equation*}
By~\cite[\S 1.8]{gilkey}, given that the fixed point set is $0$-dimensional and
that the differential operator $F$ is of degree $0$,
the asymptotic expansion of 
$\, \Tr \left(V^{-1}_\chi\, F \, e ^{-t \, \Dirac^2} \right) $ as 
 $\,  t \rightarrow 0$  involves only positive half-integer powers $\, t^{\frac{n}{2}} $, 
with $n \geq 0$. Thus, the residue at $s=0$ of the analytic continuation of the
above expression is equal to $0$.
 \smallskip 

It follows, firstly, that the cocycle $\Psi_1$ is localized
to the identity. Secondly, for $\, \psi = \phi^{-1}$,
 one obtains  the following explicit formula: 
\begin{equation*} \label{psi1}
\Psi_1 (f\,U^*_\phi , \, g\, U^*_\psi )  \, = \,
 \frac{2}{i}  \int_{\R/\Z} \, f (x)  \, \phi'(x)^{\frac{1}{2}}\,
 \big(g\, (\psi')^{-\frac{1}{2}} \big)^\prime (\phi (x)) \, dx \, .
 \end{equation*}
An elementary manipulation of this expression gives
\begin{equation} \label{psi2}
 \Psi_1 (f\,U^*_\phi , \, g\, U_\phi ) \,  =\,
  -2i \,  \int_{\R/\Z}  \, f \, \phi^\ast (dg) \,  -\, i \,  \int_{\R/\Z}  \, 
f  \, \phi^\ast (g) \, d \log \phi'   \, .
 \end{equation}
Both terms in the right hand side are recognizable cyclic cocycles, 
\cf \cite{transverse} and \cite[III. 6.$\beta$]{book}.
The first equals $\, -2i$ times the cyclic cocycle
\begin{equation} \label{tfc1}
\tau (f\,U^*_\phi , \, g\, U^*_\psi )\, = \, \varphi (f\,U^*_\phi \cdot dg\, U^*_\psi ) \, ,
 \end{equation}
with $\varphi$ the canonical state \eqref{state}, which
represents the \emph{transverse fundamental class} of the `quotient space'
$[S^1/\Gamma]$. The second coincides with    
the \emph{Lie derivative} of $\, \tau$ with respect to the generator $\delta$ of
the modular automorphism group $\{\sigma_t \}_{t \in \R}$ (see  \eqref{dfesigma1t}),
$\delta (f\,U^*_\phi) \, = \, i \, \log \phi' \cdot f\,U^*_\phi$.
Indeed, for $\, \psi = \phi^{-1}$, 
\begin{eqnarray*} \label{Lder}
\cL_\delta \, \tau (f\,U^*_\phi , \, g\, U^*_\psi )&:=& \, 
\tau (\delta(f\,U^*_\phi) , \, g\, U^*_\psi ) \, + \, \tau (f\,U^*_\phi , \, \delta (g\, U^*_\psi) ) \\ 
&=& \varphi (\delta(f\,U^*_\phi) \cdot dg\, U^*_\psi ) \, + \, 
\varphi (f\,U^*_\phi \cdot d(\delta (g\, U^*_\psi )) \\
&=& i  \int_{\R/\Z}  \, f \, \log \phi'   \, \phi^\ast (dg)  \, +
 \, i  \int_{\R/\Z}  \, f  \, \phi^\ast (g\, d \log \psi') \\
&=&- i \,  \int_{\R/\Z}  \, 
f  \, \phi^\ast (g) \, d \log \phi'   \, .
 \end{eqnarray*}
By the cyclic analogue of the Cartan homotopy formula 
$\, \cL_\delta \, = \, [e_\delta  + E_\delta , b + B]$,
and since $\tau$ is a cyclic cocycle, one sees that $ \cL_\delta \, \tau$ is a 
coboundary:
\begin{eqnarray} \label{Lder2}
 \cL_\delta \, \tau  \, &=& \, B (e_\delta \, \tau ) + b (E_\delta \, \tau)  \, , \\ \nonumber
\text{where}  \quad
e_\delta \tau (a^0 , a^1, a^2) &:=& - \tau (\delta (a^2) a^0, a^1) \quad
\text{and }  \quad E_\delta  \tau (a^0):=  \tau (1, a^0)  =  B  \tau (a^0) = 0 .
 \end{eqnarray}
We have thus proved the following statement, showing that the
 local Hochschild cocycle \eqref{oddtHcocycle} does represent the `correct' class.
 \medskip

 \begin{thm} \label{tHclass1}
 The local Hochschild cocycle associated to the 
 Dirac spectral triple over  $\cA=\,C^\infty(S^1)\rtimes \Gamma$,
 with $\Gamma$ acting by orientation preserving diffeomorphisms
 with nondegenerate isolated fixed points, 
 has the explicit expression
  \begin{equation} \label{psi3}
 \Psi_1  \,  =\, -2i \, \tau  \, + \,  \cL_\delta \, \tau  .
 \end{equation}
In particular, it is a cyclic cocycle and its
 periodic cyclic cohomology  class is proportional to the
 transverse fundamental class $[S^1/\Gamma]$.
 \end{thm}
\medskip

The facility of the above computation stands in sharp contrast to the tremendous
amount of calculations required to explicitly compute the local index
cocycle for the type II -- but $3$-dimensional -- lift of the above spectral
 triple to the frame bundle. Those calculations were kept unpublished,
due to their excessive length, but we recorded their end result in \cite[Appendix, eq. (9)]{C-M3};
 the expression of the corresponding
$3$-dimensional cocycle is quite similar to equation \eqref{psi3}. 
\medskip

\section{Future developments}

We conclude by listing, roughly in their increasing order of complexity, a few themes
for future research in this direction.
\smallskip

\subsection{Symbolic calculus and local index formula} The symbolic calculus
developed for spectral triples (\cf~\cite[Appendix B]{C-M2}) needs to be adapted to allow,
in particular, establishing that if
$(\cA , \fH , D)$ is a $\sigma$-spectral triple with
$D^{-1}\in \cL^{n,\infty}$ that satisfies a stronger regularity assumption, then
\begin{equation}  \label{extrareg}
 |D|^{-t} \left(|D|^t \, a \, - \, \sigma^t (a) \, |D|^t \right)
\in \cL^{n,\infty} \, , \qquad \forall \, t \in \R \, .
\end{equation}
By Remark \ref{orderdecrease}, the transverse spectral triple example fulfills this
property. Note also that the extra regularity assumption \eqref{extrareg}
immediately reconciles the two definitions given above to
the Chern character, {\it viz.} \eqref{CCh} and \eqref{Cycocycle}. Indeed, one 
can then produce the following homotopy between the cyclic cocycles
$\Phi_F$ and $\Phi_{D, \sigma}$:
 \begin{equation*}
 \Phi_t (a^0 , a^1, \ldots , a^n):= 
 \Tr \,\big(\gamma\,  D_t^{-1} (D_t \, a^0 - \sigma^{1-t}(a^0) \, D_t)\,
 \cdots \, D_t^{-1} (D_t \, a^n - \sigma^{1-t}(a^n) \, D_t)\big),
 \end{equation*}
 where $\, D_t = D\, |D|^{-t}$ and $\, t \in [0, 1]$.
 \smallskip

The full expression of the local formula for the Chern character of a finitely
summable $\sigma$-spectral triple,  based on $\sigma$-twisted commutators
and extending the noncommutative local index formula in~\cite[Part II]{C-M2},
remains to be worked out.  
 
\medskip

\subsection{Relation to type II and Thom isomorphism}
One should expect that the constructions in the
foliation context and in the context of modular forms of
hypoelliptic spectral triples on frame bundles extend to the general
context of twisted spectral triples satisfying (1PG). The
noncommutative space associated to the total space of the frame
bundle corresponds to the cross product algebra by the one-parameter
group $ \{ \sigma_t \}_{t \in \R} $. The $K$-homology classes on the base and on the
total space as well as their local index cyclic cocycles should be
related by a Thom isomorphism.
\medskip

\subsection{Higher dimensions} The most challenging task ahead consists
in extending the above considerations to the case of
foliations of higher codimension, which has been put in the
framework of a higher form of Tomita's theory in \cite{transverse},
Section 3. The above notion of twisting only allows to handle the
determinant part of the cocycle given by the Jacobian. One expects
the general case to involve dual actions of Lie groups such as
$\GL(n)$ and more generally of quantum groups.
\medskip

\subsection{Relation with quantum groups}
The domain of quantum groups is a natural arena where twisting
frequently occurs (see~\cite{krahmer}) and where the above extension
of the notion of spectral triple could be useful. One would expect
that the higher dimensional generalizations alluded to above
 would also extend to the braided context that arises
from quantum groups.


\bigskip

\begin{thebibliography}{999}
 
 \bibitem{bost} {\sc J.-B. Bost}, Principe d'Oka, $K$-th\'eorie et syst\`emes dynamiques
non commutatifs,
{\it Invent. math.} {\bf 101} (1990), 261-333.

\bibitem{bourg} {\sc J.-P. Bourguignon, P. Goduchon}, Spineurs, op\'erateurs de Dirac
et variations,
{\it Commun. Math. Phys.} {\bf 144} (1992), 581-599.

\bibitem{ncdg} {\sc A. Connes}, Noncommutative differential geometry,
\textit{Inst. Hautes Etudes Sci. Publ. Math.} \textbf{62} (1985),
257-360.

\bibitem{transverse} {\sc A. Connes}, Cyclic cohomology and the
transverse fundamental class of a foliation. In {\bf Geometric methods
in operator algebras} (Kyoto, 1983), pp. 52-144,
{\it  Pitman Res. Notes in Math., {\bf 123}, Longman, Harlow}, 1986.
 \url{ftp://ftp.alainconnes.org/transfund.pdf}

\bibitem{book} {\sc A. Connes}, {\bf Noncommutative geometry},
 Academic Press,1994.
 \url{ftp://ftp.alainconnes.org/book94bigpdf.pdf}
 
\bibitem{C-M2} {\sc A. Connes, H. Moscovici},
 The local index formula in noncommutative
geometry, {\it Geom. Funct. Anal.} {\bf 5} (1995), 174-243.

\bibitem{C-M3} {\sc A. Connes, H. Moscovici}, Hopf algebras, cyclic 
cohomology and the transverse index theorem, \textit{Commun. Math. 
Phys.} \textbf{198} (1998), 199-246 . 


\bibitem{C-M} {\sc A. Connes, H. Moscovici}, Cyclic cohomology and Hopf algebras,
{\it Lett. Math. Phys.} {\bf 48} (1999), 97-108.

 \url{ftp://ftp.alainconnes.org/letters.pdf}

\bibitem{gilkey} {\sc P. B. Gilkey},
{\bf Invariance theory, the heat equation, and the Atiyah-Singer index theorem},
 Mathematics Lecture Series, 11. Publish or Perish, Inc., Wilmington, DE, 1984.

\bibitem{krahmer} {\sc  Tom
Hadfield, Ulrich Kraehmer},  On the Hochschild homology of quantum
SL(N),  math.QA/0509254.

\bibitem{hitchin} {\sc N. Hitchin}, Harmonic spinors,
{\it Adv. Math.} {\bf 14} (1974), 1-55.


\end{thebibliography}
\end{document}